\documentclass[a4paper,10pt,reqno]{amsart}

\usepackage{amsthm,amsmath,amsfonts,amssymb}
\usepackage{pgf, tikz-cd}
\usepackage{mathrsfs}
\usetikzlibrary{arrows}
\definecolor{cqcqcq}{rgb}{0.7529411764705882,0.7529411764705882,0.7529411764705882}
\def\step#1#2{\par\noindent{\underline{\it Step~#1.}}\emph{ #2}\\}

\newcounter{mt}
\def\maintheorem#1#2#3{\par \medskip \noindent {\bf Theorem~\mref{#1}}~(#2).~{\it #3}\par}
\def\mref#1{\Alph{#1}}
\def\maintheoremdeclaration#1{\stepcounter{mt}\newcounter{#1}\setcounter{#1}{\arabic{mt}}}

\maintheoremdeclaration{Main}

\def\XXint#1#2#3{{\setbox0=\hbox{$#1{#2#3}{\int}$} \vcenter{\vspace{-1pt}\hbox{$#2#3$}}\kern-.5\wd0}}
\def\Xint#1{\mathchoice {\XXint\displaystyle\textstyle{#1}}{\XXint\textstyle\scriptstyle{#1}}{\XXint\scriptstyle\scriptscriptstyle{#1}}{\XXint\scriptscriptstyle\scriptscriptstyle{#1}}\!\int}
\def\intmed{\hbox{\ }\Xint{\hbox{\vrule height -0pt width 9pt depth 0.7pt}}}

\numberwithin{equation}{section}

\def\eps{\varepsilon}
\def\C{\mathcal{C}}
\def\R{\mathbb{R}}
\def\N{\mathbb{N}}
\def\S{\mathbb{S}}
\def\SS{\mathcal{S}}
\def\H{\mathcal{H}}
\def\Lloc{ L^1_{\text{loc}}(\R^n)}
\def\J{\mathcal{J}}

\theoremstyle{plain}
\newtheorem{thm}{Theorem}[section]
\newtheorem{lem}[thm]{Lemma}
\newtheorem{prop}[thm]{Proposition}

\theoremstyle{definition}
\newtheorem{defin}[thm]{Definition}

\title[]{On the isoperimetric problem with double density}

\author{Aldo Pratelli}
\address{Department Mathematik, Universit\"at Erlangen-N\"urnberg, Cauerst. 11, 91058 Erlangen - Germany}
\email{pratelli@math.fau.de}
\author{Giorgio Saracco}
\address{Department Mathematik, Universit\"at Erlangen-N\"urnberg, Cauerst. 11, 91058 Erlangen - Germany}
\email{saracco@math.fau.de}

\thanks{G. Saracco has been supported by the DFG Grant no. GZ:PR 1687/1-1}

\begin{document}


\begin{abstract}
In this paper we consider the isoperimetric problem with double density in Euclidean space; that is, we study the minimisation of the perimeter among subsets of $\R^n$ with fixed volume, where volume and perimeter are relative to two different densities. The case of a single density, or equivalently, when the two densities coincide, has been well studied in the last years; nonetheless, the problem with two different densities is an important generalisation, also in view of applications. We will prove the existence of isoperimetric sets in this context, extending the known results for the case of single density.
\end{abstract}

 \hspace{-3cm}
 {
 \begin{minipage}[t]{0.6\linewidth}
 \begin{scriptsize}
 \vspace{-2.5cm}
 This is a pre-print of an article published in \emph{Nonlinear Analysis} (2018). The final authenticated version is available online at: https://doi.org/10.1016/j.na.2018.04.009
 \end{scriptsize}
\end{minipage} 
}


\maketitle


\section{Introduction}

For any $n\geq 2$, we consider the isoperimetric problem in $\R^n$ with double density. That is, two l.s.c. and locally summable functions $f:\R^n\to \R^+$ and $h:\R^n\times \S^{n-1}\to\R^+$ are given, and we measure the volume and perimeter of any Borel set $E\subseteq\R^n$ according to the formulas
\begin{align*}
|E|_f := \int_E f(x)\, dx\,, && P_h(E):= \int_{\partial^* E} h(x,\nu(x))\, d\H^{n-1}(x)\,,
\end{align*}
where as usual $\partial^* E$ is the reduced boundary of $E$ and for every $x\in\partial^* E$ the vector $\nu(x)\in\S^{n-1}$ is the outer normal at $x$, see for instance~\cite{AFP}. Of course, the standard Euclidean case corresponds to the situation where $f$ and $h$ are constantly equal to $1$. We will refer to the ``single density'' case when $h(x,\nu)=f(x)$ for every $x\in\R^n,\,\nu\in\S^{n-1}$.\par

In recent years, for different reasons, there was much attention on the study of problems in $\R^n$ endowed with a single density, a non-exhaustive list of references is~\cite{BBCLT, CMV,CRoS,DFP,DHHT,MP,RCBM}. Nevertheless, as already pointed out by several authors, the case of the double density is an important generalisation, since many of the possible applications correspond to two different densities. The simplest example is given by the Riemannian manifolds: of course they locally behave like $\R^n$ with double density, being $f$ the norm of the Riemannian metrics, and $h$ their derivatives. In particular, it is important that the density $h$ relative to the perimeter depends not only on a point, but also on the direction of the tangent space at that point. A preliminary study on an isoperimetric-like problem with double density is for instance the paper~\cite{SarMAMA}.

When we refer to the ``isoperimetric problem'', as always we mean that we are interested in minimising the (relative) perimeter of sets with fixed (relative) volume. Hence, we define as usual the \emph{isoperimetric profile} as
\begin{equation}\label{eq:iso1}
\J(V) := \inf \left \{ P_h(E)\,:\, |E|_f=V\right\}\,,
\end{equation}
and we are interested in studying the existence of isoperimetric sets $E$, that is, the sets realising the infimum above. Let us immediately notice that the existence problem becomes interesting only when $f\notin L^1(\R^n)$. Indeed, given any volume $V>0$ we can always take a minimising sequence $\{E_k\}_k$ relative to the volume $V$; up to subsequences, the sequence locally converges in $L^1$ to some set $E$, and by the lower semicontinuity of $h$ we immediately have that $P_h(E)\leq \liminf_k P_h(E_k)=\J(V)$. The set $E$ is then obviously an isoperimetric set of volume $V$ as soon as its volume is $V$; when $f\in L^1(\R^n)$ this is always true, so the existence issue is always trivial. Instead, if $f\notin L^1(\R^n)$, it is only true that $|E|_f\leq V$, but some mass could have ``escaped at infinity''. It is easy to see that, in general, the existence of isoperimetric sets fails, and it is interesting to understand under which conditions this is true.\par
For the case of the single density much is known. More precisely, since as noticed above there is no existence issue when the volume of $\R^n$ is finite, the main point is to study the behaviour of not $L^1$ densities at infinity. Roughly speaking in order to minimise the perimeter, the sets ``prefer'' the zones where the density is low. It is then reasonable to expect that minimising sequences remain bounded if the density explodes at infinity, which leads to no loss of volume, hence to existence; on the other hand, one expects that such sequences escape at infinity if the density goes to $0$ at infinity, which leads to non-existence as $\J$ would be identically $0$ and no set has zero perimeter, at least for strictly positive densities. In fact, in~\cite{MP} it was proved that existence is true for the case of radial densities which explode at infinity (the radial assumption is not needed if $n=2$, but necessary if $n\geq 3$ as counterexamples show). Moreover, existence generally fails for densities going to $0$ at infinity, unless they belong to $L^1(\R^n)$, as previously noticed. Hence, the case when the (single) density goes to $0$ or explodes at infinity is more or less known. Additionally, simple examples show that existence may fail if the $\liminf$ and the $\limsup$ of the density at infinity are different. Then, the only interesting case which remains to describe is when the density goes to a non-zero limit at infinity. The final answer about this point was given in~\cite{DFP}: there, the existence of isoperimetric sets of any volume was proved under the assumption that the density converges to its limit from below in the sense of Definition~\ref{frombelow}, see Theorem~\ref{DFP}. On the other hand, easy examples show that existence is generally false if the density does not converge from below, hence we can say that the interesting issues about existence in the single density case are all solved.
\begin{defin}\label{frombelow}
Let $\sigma:\R^n\to\R^+$ be a l.s.c. function. We say that \emph{$\sigma$ converges from below to $\ell>0$} if $\sigma(x)\to \ell$ when $|x|\to \infty$, and there exists $M\in\R$ such that $\sigma(x)\leq \ell$ whenever $|x|>M$. Analogously, we say that \emph{$\sigma$ converges from above to $\ell>0$} if $\sigma(x)\to\ell$ when $|x|\to\infty$ and there exists $M\in\R$ such that $\sigma(x)\geq \ell$ whenever $|x|>M$.
\end{defin}
\begin{thm}\label{DFP}
Let $f:\R^n\to \R^+$ be a l.s.c., locally integrable function, converging from below to a non-zero limit. Then, for every positive volume $V>0$ there exist isoperimetric sets, that is, Borel sets $\overline E\subseteq\R^n$ which minimise $P_f(E)$ among all the sets satisfying $|E|_f=V$.
\end{thm}

Let us now pass to discuss the general case of double densities. Arguing similarly as before, one can see that the most interesting case is when both $f$ and $h$ converge to a non-zero limit at infinity. Up to a reparametrisation and for ease of notations, we will assume that they both converge to $1$. Let us then define the functions $h^+,\,\tilde f,\,\tilde h:\R^n\to\R^+$ as
\begin{align}\label{deft}
h^+(x)=\sup_{\nu\in\S^{n-1}} h(x,\nu)\,, && \tilde f(x)=|f(x)-1|\,, && \tilde h(x)=|h^+(x)-1|\,,
\end{align}
and notice that $\tilde f$ and $\tilde h$ converge to zero at infinity. It is easy to see that isoperimetric sets always exist if $f$ converges to $1$ from above and $h^+$ from below, while existence generally fails if $f$ converges from below and $h^+$ from above. We are thus led to consider the case when both $f$ and $h^+$ converge from below, or they both converge from above. Our main result is the following.
\maintheorem{Main}{Isoperimetric sets in $\R^n$ with double density}{Let $f:\R^n\to\R^+$ and $h:\R^n\times\S^{n-1}\to\R^+$ be two l.s.c. and locally integrable functions, both converging to $1$ at infinity. There exist isoperimetric sets of any volume in the following two cases:
\begin{itemize}
\item both $f$ and $h^+$ converge from below to $1$, and
\begin{equation}\label{eq:HPbelow}
\inf \big\{C\geq 0:\, \tilde f(x)\leq C \tilde h(x) \hbox{ for } |x|\gg 1\big\} < \frac{n}{n-1};
\end{equation}
\item both $f$ and $h^+$ converge from above to $1$,
\begin{equation}\label{eq:HPabove}
\sup \big\{C\geq 0:\, \tilde f(x)\geq C \tilde h(x) \hbox{ for } |x|\gg 1\big\} > \frac{n}{n-1}\,,
\end{equation}
and for any $R\gg 1$ one has $\int_R^\infty \tilde h_r(t)\,dt=+\infty$, where the subscript $r$ denotes the radial average, i.e.
\begin{equation}\label{defra}
\tilde h_r(t) := \intmed_{\partial B_{|t|}} \tilde h(x)\,d\H^{n-1}(x)\,.
\end{equation}
\end{itemize}}

A couple of comments are in order. First of all, as already pointed out above, the very same result is true if $f$ and $h$ converge to two limits $a,\,b>0$; in fact, one can always reduce to the case $a=b=1$ up to a multiplication of $f$ and $h$ by a constant, which obviously does not change isoperimetric sets. Moreover, notice that Theorem~\mref{Main} generalises Theorem~\ref{DFP}: in fact, in the case of the single density one has $\tilde f/\tilde h=1<n/(n-1)$, so according with the result of Theorem~\ref{DFP} one gets existence for the case of convergence from below, and not for the case of convergence from above. It is then interesting that existence is obtained also in some cases when the convergence is from above, unlike the case of the single density. Observe also that the convergence to $1$ from below or from above is required for $h^+$, not for $h$. However, if $h^+$ converges from below to $1$, then in particular $h(x,\nu)\leq 1$ for every $x$ with $|x|\gg 1$ and every $\nu\in\S^{n-1}$; on the other hand, if $h^+$ converges from above to $1$, then $h(x,\nu)<1$ is possible for infinitely large $|x|$ and some directions $\nu\in\S^{n-1}$. Finally, notice the additional assumption about the  non-summability of $\tilde h_r$ in the case of convergence from above: surprisingly enough, this assumption is really needed, as our example in Section~\ref{sec:counter} shows.\par\smallskip

The plan of the paper is very simple; in Section~\ref{sec:prelim} we collect a couple of technical properties, which are well-known for the case of single densities and that we generalise to the case of double density. Then, in Section~\ref{sec:proof} we present the proof of Theorem~\mref{Main}, which is a careful extension of the argument already used in~\cite{DFP} to prove Theorem~\ref{DFP}. Finally, in Section~\ref{sec:counter} we show an example in which existence fails for the case of two radial densities $f,\,h=h^+$ which converge from above to $1$ and for which~(\ref{eq:HPabove}) holds true, but for which $\tilde h_r$ is summable.

\section{Preliminary properties\label{sec:prelim}}

In this section, we present the definition of ``mean density'' and a couple of basic results: both the definition and the results are well-known for the case of a single density (for the definition, see for instance~\cite{MP}, for the results, see for instance~\cite[Lemma~2.1,2.3]{DFP}); here, we generalise them to the case of double density. In the following, $f:\R^n\to\R^+$ and $h:\R^n\times \S^{n-1}\to\R^+$ always denote a double density, that is, two l.s.c. and locally summable functions. Moreover, for the sake of brevity, whenever a function $g:\R^n\times\S^{n-1}\to\R^+$ and a $(n-1)$-dimensional oriented surface $\Gamma$ are given, we write
\[
\H^{n-1}_g(\Gamma) = \int_{x\in\Gamma} g(x,\nu(x))\,d\H^{n-1}(x)\,,
\]
where $\nu(x)\in\S^{n-1}$ is the outer (in the sense of the orientation of $\Gamma$) normal vector to $\Gamma$ at $x$.

\begin{defin}\label{meandensity}
Let $F\subseteq\R^n$ be a set of finite perimeter. We define the $(f,h)$-mean density of $F$ as the number $\rho\in\R^+$ such that
\[
P_h(F) = n(\omega_n \rho)^\frac 1n |F|_f^{\frac{n-1}{n}}\,.
\]
\end{defin}

It is important to realise the meaning of the ``mean density''. Namely, the $(f,h)$-mean density of $F$ is the unique number $a\in\R^+$ such that, for the case of $f\equiv h\equiv a$, a ball with volume $|F|_f$ has precisely perimeter $P_f(F)$.

\begin{lem}\label{lem:isoE}
Suppose that, far enough from the origin, $h$ is locally bounded from above and $h^+\leq \lambda f$ for some $\lambda>0$. Let $V>0$ be given, and let $\{E_j\}_j$ be an isoperimetric sequence corresponding to volume $V$, which converges in $\Lloc$ to $E \subseteq\R^n$. Then, $E$ is isoperimetric for the volume $|E|_f$. Moreover, if $f$ and $h$ pointwise converge to some $a>0$ at infinity, then
\begin{equation}\label{eq:scarto}
\J(V) = P_h(E) +n(a\omega_n)^{\frac1n}(V-|E|_f)^{\frac{n-1}{n}}\,.
\end{equation}
\end{lem}
\begin{proof}
Notice that, since $f$ and $h$ are l.s.c., then
\begin{align}\label{eq:liminfperimetro}
|E|_f\leq V\,, && P_h(E) \le \liminf_j P_h(E_j) = \J(V)\,.
\end{align}
We divide the proof in two parts.
\step{1}{The set $E$ is isoperimetric for its volume.}
We have to prove that $P_h(E)=\J(|E|_f)$. This is emptily true if $|E|_f=0$, and it is given by~(\ref{eq:liminfperimetro}) if $|E|_f=V$. We are then left to consider the case $0<|E|_f<V$, and we assume by contradiction the existence of a set $F$ with $|F|_f = |E|_f$ and $P_h(F)<P_h(E)$. For simplicity of notation, we define
\begin{equation}\label{eq:eta}
\eta:=\frac{P_h(E)-P_h(F)}5\,.
\end{equation}
Choose $x\in\R^n$ a Lebesgue point for the set $F$ and for both the functions $f$ and $h$: such a point exists as $|F|_f=|E|_f>0$. As a consequence, we can find $r\ll 1$ such that
\begin{align}\label{eq:lebfh}
\frac12\,\omega_n f(x) r^n \leq |B_r(x) \cap F|_f \leq 2\omega_n f(x) r^n\,,&& P_h(\R^n\setminus B_r(x)) \leq 2n\omega_n h^+(x) r^{n-1}\,.
\end{align}
In particular, the first estimate is true for every $r\ll 1$, and the second for almost every $r\ll 1$.\par

Similarly, let $y\in\R^n$ be a Lebesgue point for the set $\R^n\setminus F$ and the functions $f$ and $h$; observe that such a point surely exists because the fact that $|F|_f<|E|_f$ implies in particular $f\notin L^1(\R^n)$; moreover, $y$ can be taken arbitrarily far from the origin, so that in particular, for some $M>0$, one has $h<M$ in a neighbourhood of $y$. In particular, there exists $\bar\rho\ll |y-x|$ such that, for every $\rho<\bar\rho$, one has
\begin{align}\label{eq:aim2}
|B_\rho(y) \setminus F|_f \geq \frac 12\,\omega_n f(y)\rho^n\,, && P_h(B_\rho(y)) \leq Mn\omega_n\rho^{n-1}<\eta\,.
\end{align}
Let us now take $\eps<\eta$ in such a way that
\begin{equation}\label{etasmall}
\omega_n f(y) \bar \rho^n > 2\eps\,.
\end{equation}
We claim the existence of a set $F'$ and of a big radius $R>|y|+\bar\rho$ such that
\begin{align}\label{eq:aim}
F'\cap B_{\bar\rho}(y)=F\cap B_{\bar\rho}(y)\,, && F'\subseteq B_R\,, && P_h(F')<P_h(E)-4\eta\,, && 0<\eps' := |E|_f -|F'|_f < \eps\,.
\end{align}
In fact, if $F$ is bounded then it is enough to define $F'=F\setminus B_r(x)$ for a suitably small $r$, keeping in mind~(\ref{eq:lebfh}) and the fact that $P_h(F\setminus B_r(x))\leq P_h(F)+P_h(\R^n\setminus B_r(x))$. Otherwise, if $F$ is not bounded, we set $F'=F\cap B_R$ for some $R\gg 1$. The set equality and the inclusion in~(\ref{eq:aim}) are true by construction, as well as the inequalities about the volume if $R\geq R_1$ for some $R_1$ big enough. So, we have to prove the inequality about the perimeters. Since $R_1$ is an arbitrarily large constant, we can assume that $h^+\leq \lambda f$ on $\R^n\setminus B_{R_1}$. Suppose then for a moment that it is false for every $R\geq R_1$, so that
\[
P_h(F\cap B_R)>P_h(F)+\eta\,.
\]
Then, since for almost every $R\geq R_1$ we have
\[
P_h(F\cap B_R) = P_h(F) + \int_{F\cap\partial B_R} h(x,x/|x|)\, d\H^{n-1}(x) - \int_{\partial F\setminus B_R} h(x,\nu(x))\, d\H^{n-1}(x)\,,
\]
and we deduce
\[\begin{split}
V&>|F|_f > |F\setminus B_{R_1}|_f
= \int_{R_1}^\infty \int_{F\cap \partial B_R} f(x)\, d\H^{n-1}(x) \, dR
\geq \frac 1\lambda\, \int_{R_1}^\infty \int_{F\cap \partial B_R} h(x,x/|x|)\,d\H^{n-1}(x)\,dR\\
&\geq \frac 1\lambda\, \int_{R_1}^\infty \eta\,dR= \infty\,,
\end{split}\]
which is a contradiction. The existence of a constant $R$ and a set $F'$ satisfying~(\ref{eq:aim}) is then established.\par
Let us now consider the set $E$. It is clearly possible to fix some $R'>R$ such that
\begin{align}\label{toadd}
|E\setminus B_{R'}|_f < \frac{\eps'}{2\lambda} \,,&& \int_{\partial^*E \setminus B_{R'}} h(x,\nu(x))\,d\H^{n-1} < \eta\,.
\end{align}
Keep in mind that, since $E_j$ converges in $\Lloc$ to $E$, then $E_j \cap B_{R'}$ and $E_j \cap B_{R'+1}$ converge in $L^1$ to $E \cap B_{R'}$ and $E \cap B_{R'+1}$ respectively. Therefore, by~(\ref{toadd}), recalling also the lower semicontinuity of $h$, for every sufficiently large index $j$ we have
\begin{gather}
|E|_f - \frac {\eps'}{1+2\lambda} < |E_j \cap B_{R'}|_f \le |E_j \cap B_{R'+1}|_f < |E|_f +\frac{\eps'}{1+2 \lambda}\,,\label{eq:volumi}\\
P_h(E) < \int_{\partial^* E_j\cap B_{R'}} h(x,\nu_j(x))\,d\H^{n-1} + \eta\,, \label{eq:perimetri}
\end{gather}
where by $\nu_j(x)$ we denote the outer direction of the boundary $\partial^* E_j$ at $x\in\partial^* E_j$. By~\eqref{eq:volumi}, we get
\[
\frac\eps\lambda > \frac{\eps'}\lambda \geq |E_j \cap(B_{R'+1}\setminus B_{R'})|_f =\int_{R'}^{R'+1}\H^{n-1}_f (E_j\cap\partial B_R)\, dR\,,
\]
from which we deduce the existence of some radius $R_j \in (R', R'+1)$ such that
\begin{equation}\label{eq:stima4}
\H^{n-1}_{h^+} (E_j \cap \partial B_{R_j}) \leq \lambda \H^{n-1}_f (E_j \cap \partial B_{R_j}) < \eps\,.
\end{equation}
Next, we define the set $G_j = F' \cup (E_j \setminus B_{R_j})$; keeping in mind that $R_j>R'>R$, (\ref{eq:aim}), (\ref{eq:volumi}) and the fact that $|E_j|_f=|E|_f=V$ by definition, we estimate the volume of $G_j$ by
\begin{equation}\label{eq:volGj}
|G_j|_f = |F'|_f +|E_j \setminus B_{R_j}|_f = |E|_f -\eps' +|E_j \setminus B_{R_j}|_f \in (V-\eps, V)\,,
\end{equation}
and its perimeter, thanks to~(\ref{eq:aim}), (\ref{eq:stima4}) and~(\ref{eq:perimetri}) and the fact that $R_j>R'$, by
\[\begin{split}
P_h(G_j) &= P_h(F') +P_h(E_j \setminus B_{R_j}) 
< P_h(E) -4\eta + \int_{\partial^*E_j \setminus B_{R_j}} h(x,\nu_j(x))\,d\H^{n-1} + \H_{h^+}^{n-1}(E_j\cap \partial B_{R_j})\\
&< P_h(E) -4\eta + P_h(E_j)-\int_{\partial^*E_j \cap B_{R'}} h(x,\nu_j(x))\,d\H^{n-1} + \eps
<  P_h(E_j) -3\eta+ \eps<P_h(E_j) -2\eta\,.
\end{split}\]
Finally, let us take $\rho_j<\bar\rho$ in such a way that the set $D_j = G_j \cup B_{\rho_j}(y)$ satisfies $|D_j|_f=V$, which is possible by~(\ref{eq:volGj}), (\ref{eq:aim2}), (\ref{eq:aim}) and~(\ref{etasmall}). Putting together the last estimate and again~(\ref{eq:aim2}), we get
\[
P_h(D_j) \leq P_h(G_j) + P_h(B_{\rho_j}(y))
< P_h(E_j) -\eta\,.
\]
Since $\{E_j\}$ is a minimising sequence for $\J(V)$, we deduce that $P_h(D_j)<\J(V)$ for $j\gg 1$, which is a contradiction. This concludes the first part of the claim.

\step{2}{The validity of~(\ref{eq:scarto}).}
We now assume that $f$ and $h$ converge to some $a>0$ at infinity, and we have to prove the validity of~\eqref{eq:scarto}. We can assume that $|E|_f < V$, since otherwise the claim is emptily true. Arguing as in the first step, for every $\eps\ll 1$ we can find a big radius $R$ and a set $F\subseteq B_R$ such that
\begin{align*}
|F|_f \geq |E|_f -\eps\,, && P_h(F)\leq P_h(E) +\eps\,.
\end{align*}
Now, take a ball $B=B(x,r)$ having volume $|B|_f=V- |F|_f$ and very far from the origin. In particular, we can assume that $B\cap B_R=\emptyset$, and that the values of $f$ and $h$ in a neighbourhood of $B$ are all between $a-\eps$ and $a+\eps$. Hence,
\[
(a-\eps)\omega_nr^n \leq  |B|_f \leq V-|E|_f +\eps\,,
\]
which, since by construction $G=F\cup B$ has $f$-volume equal to $V$, yields
\[\begin{split}
\J(V) &\le P_h(G) = P_h(F) +P_h(B) \le P_h(E)+\eps +n\omega_n (a+\eps) r^{n-1}\\
&\leq P_h(E) +\eps +\frac{n\omega_n^{\frac 1n}(a+\eps)}{(a-\eps)^{\frac{n-1}{n}}}  \left( V-|E|_f+\eps \right)^{\frac {n-1}{n}}\,.
\end{split}\]
By letting $\eps\to 0$ we obtain the first inequality in~\eqref{eq:scarto}.\par

To show the opposite inequality, for every $\eps>0$ we can argue as usual to pick a large $R$ so that
\begin{align*}
|E\cap B_R|_f > |E|_f -\eps\,, && P_h(E\setminus B_R)< \eps\,.
\end{align*}
As in the first step, for every $j\gg 1$ we can find some $R_j \in (R, R+1)$ such that
\begin{align*}
|E_j \cap B_{R_j}|_f < |E|_f +\eps\,, && \H^{n-1}_{h^+} (E_j \cap \partial B_{R_j}) < \eps\,, && P_h(E) <P_h(E_j \cap B_{R_j}) +\eps\,,
\end{align*}
keeping in mind that both $f$ and $h$ are arbitrarily close to $a$ outside of $B_R$, up to taking a sufficiently large $R$. Then, calling $P(\cdot)$ and $|\cdot|$ the standard Euclidean perimeter and volume, by the standard isoperimetric inequality we can estimate
\[\begin{split}
P_h(E_j \setminus B_{R_j}) &> (a-\eps)P(E_j \setminus B_{R_j}) \ge (a-\eps) n\omega_n^{\frac 1n}|E_j\setminus B_{R_j}|^{\frac{n-1}{n}}\\
&> \frac{(a-\eps)n\omega_n^{\frac 1n}}{(a+\eps)^{\frac{n-1}{n}}} \, |E_j \setminus B_{R_j}|_f^{\frac{n-1}{n}}
>\frac{(a-\eps)n\omega_n^{\frac 1n}}{(a+\eps)^{\frac{n-1}{n}}} \, \big(V-|E|_f -\eps \big)^{\frac{n-1}{n}}\,.
\end{split}\]
As a consequence, we deduce
\begin{align*}
P_h(E_j) &\geq P_h(E_j \cap B_{R_j}) +P_h(E_j \setminus B_{R_j}) -2\H_{h^+}^{n-1}(E_j \cap \partial B_{R_j}) \\
&\ge P_h(E) -3\eps + \frac{(a-\eps)n\omega_n^{\frac 1n}}{(a+\eps)^{\frac{n-1}{n}}}  \,\big(V-|E|_f -\eps \big)^{\frac{n-1}{n}}\,.
\end{align*}
Since $E_j$ is an isoperimetric sequence for volume $V$, thus $P_h(E_j)\to \J(V)$, by letting $\eps \to 0$ we get the second inequality in~(\ref{eq:scarto}), so the thesis is concluded.
\end{proof}

\begin{lem}\label{lem:boundedE}
Assume that the densities $f$ and $h$ pointwise converge to $a>0$ at infinity, let $\{E_j\}_j$ be a minimising sequence for some volume $V>0$ converging in $\Lloc$ to a set $E$. If $|E|_f < V$, then $E$ is bounded.
\end{lem}
\begin{proof}
For every $t>0$, we define $m(t)$ the mass of $E$ outside the ball $B_t$, that is,
\[
m(t) = |E\setminus B_t|_f = \int_t^{\infty} \H^{n-1}_f(E\cap \partial B_r)\, dr\,.
\]
Now, we pick a ball $B$ of volume equal to $V-|E|_f+m(t)$ arbitrarily far from the origin. Then, the set $E\cap B_t \cup B$ has volume exactly $V$, thus $\J(V)\leq P_h(E\cap B_t)+P_h(B)$. Since $f$ and $h$ are arbitrarily close to $a$ in a neighbourhood of $B$, arguing as in the second step of last lemma we deduce
\[
\J(V) \leq P_h(E\cap B_t) +n(a\omega_n)^{\frac 1n} \left( V-|E|_f +m(t)\right)^{\frac {n-1}{n}}\,.
\]
Comparing this inequality with~(\ref{eq:scarto}), and keeping in mind that $|E|_f<V$, we obtain
\[
P_h(E) \leq P_h(E\cap B_t) +c_1 m(t)
\]
for some geometrical constant $c_1$, only depending on $a,\, V,\, |E|_f$ and $n$. Let now $\eps\ll 1$ be fixed; as soon as $t$ is large enough, using again that $f$ and $h$ are arbitrarily close to $a$ outside of $B_t$, applying the standard isoperimetric inequality to $E\setminus B_t$, and keeping in mind that $m'(t)=-\H^{n-1}_f(E\cap \partial B_t)$, we have then
\[\begin{split}
c_1 m(t) &\geq  P_h(E) -P_h(E\cap B_t) = P_h(E\setminus B_t) - \int_{E\cap \partial B_t} h(x,x/|x|)+h(x,-x/|x|) \,d\H^{n-1}(x) \\
&\geq \frac{(a-\eps)n\omega_n^{\frac 1n}}{(a+\eps)^{\frac{n-1}{n}}}\, |E\setminus B_t|_f^{\frac{n-1}{n}}-3 \H^{n-1}_f(E\cap\partial B_t)
\geq c_2 m(t)^{\frac{n-1}{n}} +3 m'(t)\,,
\end{split}\]
where $c_2$ is any constant smaller than $n(a\omega_n)^{1/n}$. Since by definition $m(t)\searrow 0$ for $t\to\infty$, we deduce
\[
-m'(t)\geq \frac{c_2}4\, m(t)^{\frac{n-1}n}
\]
for $t$ big enough. Since every positive solution $m$ of this inequality vanishes in a finite time, we obtain $m(t)=0$ for every $t\gg 1$, that is, the set $E$ is bounded as required.
\end{proof}

\section{Existence of isoperimetric sets\label{sec:proof}}

This section is devoted to show our main result, Theorem~\mref{Main}. The key point is to get the existence of a set $F$ arbitrarily far from the origin with $(f,h)$-mean density (see Definition~\ref{meandensity}) smaller than $1$. More precisely, we will show the following two facts.

\begin{prop}\label{prop:SetMD<1below}
Let $f$ and $h$ be densities converging from below to $1$ satisfying~\eqref{eq:HPbelow}. Then, for every $m>0$ there exists a set $F$ with volume $m$ and $(f,h)$-mean density smaller than $1$ arbitrarily far from the origin.
\end{prop}

\begin{prop}\label{prop:SetMD<1above}
Let $f$ and $h$ be densities converging to $1$ and satisfying~\eqref{eq:HPabove} such that $f$ and $h^+$ converge from above to $1$ and $\int_R^{\infty }\tilde h_r(t)\,dt = +\infty$ for every $R\gg 1$. Then, for every $m>0$ there exists a set $F$ with volume $m$ and $(f,h)$-mean density smaller than $1$ arbitrarily far from the origin.
\end{prop}

We can immediately see that our main result is a very simple consequence of the above facts.

\begin{proof}[Proof of Theorem~\mref{Main}]
Let $V>0$ be given, and let $\{E_k\}_{k\in\N}$ be a minimising sequence for the isoperimetric problem with volume $V$. Up to a subsequence, the sets $E_k$ converge to some set $E\subseteq\R^n$ in the $L^1_{\rm loc}$ sense. Thanks to Lemma~\ref{lem:isoE}, we know that $E$ is an isoperimetric set for volume $|E|_f\leq V$, and that~(\ref{eq:scarto}) holds. If $|E|_f=V$, there is nothing to prove, otherwise by Lemma~\ref{lem:boundedE} we also know that $E$ is bounded. By Proposition~\ref{prop:SetMD<1below} or~\ref{prop:SetMD<1above}, we find a set $F$ with volume $|F|_f = V - |E|_f >0$ which has $(f,h)$-mean density smaller than $1$; by Definition~\ref{meandensity}, this means
\begin{equation}\label{oh}
P_h(F) \leq n\omega_n^{\frac 1n} (V-|E|_f)^{\frac{n-1}n}\,.
\end{equation}
Since $F$ can be taken arbitrarily far from the origin and $E$ is bounded, we can assume without loss of generality that $E$ and $F$ are a strictly positive distance apart. As a consequence, the set $G=E\cup F$ has exactly volume $V$, and $P_h(G)=P_h(E)+P_h(F)$: in virtue of~(\ref{eq:scarto}) and~(\ref{oh}), this means that $G$ is an isoperimetric set for volume $V$, which concludes the proof.
\end{proof}

The plan of the section is simple. Firstly, we introduce some notation which will be useful later. Then, in Section~\ref{secbel} and~\ref{secabo} we prove respectively Proposition~\ref{prop:SetMD<1below} and~\ref{prop:SetMD<1above}.\par

Given positive numbers $R,\,\delta$ and a direction $\theta \in \S^{n-1}$, we will denote by $B^{R\theta}_\delta$ the ball of radius $\delta$ centred at $R\theta$, with $B_\delta=B^0_\delta$. We will sometimes divide the ball $B$ and its boundary $\partial B$ into the ``upper'' and ``lower'' halves $B_+$ and $\partial_+ B$, and $B_-$ and $\partial_- B$. That is, we shall fix some hyperplane $H$ passing through the origin and the point $R\theta$, call $H_\pm$ the two half-spaces having $H$ as a boundary, and then write
\begin{align}\label{eq:udlrset}
B_+ := B\cap H_+\,, &&B_-:= B\cap H_-\,, &&
\partial_+ B := \partial B\cap H_+\,, &&\partial_- B:= \partial B\cap H_-\,.
\end{align}

\subsection{The case of densities converging from below\label{secbel}}

Let us start by finding a ball, arbitrarily far from the origin, whose $(1-h)$-perimeter is controlled from below by a suitable constant times the $(1-f)$-volume.
\begin{lem}\label{prop:fhbelowstime}
Let us assume that the densities $f$ and $h$ converge from below to $1$ and satisfy~\eqref{eq:HPbelow}. Then, for every $\eps\ll 1$ there exists a ball $B$ of radius $1$ arbitrarily far from the origin such that
\begin{equation}\label{eq:pallalontanobelow}
P_{1-h}(B) \geq \big(n-1+2\eps n \big) |B|_{1-f}\,.
\end{equation}
\end{lem}
\begin{proof}
Keeping in mind~(\ref{eq:HPbelow}) and the fact that $f$ and $h$ converge to $1$ from below, we can take $\eps>0$ so small that, for $x\in\R^n$ far enough from the origin,
\[
\tilde f(x)\leq \bigg(\frac n{n-1}-8\eps n\bigg)\tilde h(x)\,, 
\]
where we recall~(\ref{deft})
\[
\tilde f(x)=1-f(x)\,, \qquad \qquad \tilde h(x)=1-h^+(x)\,.
\]
As a consequence, by~(\ref{deft}), for any such $x$ and for any $\nu\in\S^{n-1}$ we have
\begin{equation}\label{ae}
1-f(x)=\tilde f(x) \leq \bigg(\frac n{n-1}-8\eps n\bigg)\, \tilde h(x)\leq \bigg(\frac n{n-1}-8\eps n\bigg)\, \big(1-h(x,\nu)\big)\,.
\end{equation}
Using again the fact that $f\leq 1$ away from the origin, we can apply~\cite[Proposition~3.2]{DFP} to get a ball with radius $1$, arbitrarily far from the origin, such that
\[
P_{1-f}(B)\geq (n-\eps) |B|_{1-f}\,.
\]
Together with~(\ref{ae}), this readily gives
\[
P_{1-h}(B)\geq \bigg(\frac{n-1}n + 2\eps n\bigg) P_{1-f}(B)
\geq \bigg(\frac{n-1}n + 2\eps n\bigg)(n-\eps) |B|_{1-f}\,,
\]
which implies~(\ref{eq:pallalontanobelow}).
\end{proof}

We are now ready to prove Proposition~\ref{prop:SetMD<1below}.
\begin{proof}[Proof of Proposition~\ref{prop:SetMD<1below}]
First of all we notice that, up to a homothety, we can for simplicity assume that $m=\omega_n$. Let then $\eps\ll 1$ be a small constant, and apply Lemma~\ref{prop:fhbelowstime} to get a ball $B$ of radius $1$ arbitrarily far from the origin and such that~(\ref{eq:pallalontanobelow}) holds. In particular, the ball is $B=B^{R\theta}_1$ for some $R\gg n/\eps$ and $\theta\in\S^{n-1}$ and, as soon as the ball is far enough from the origin, we have
\begin{align}\label{farenough}
1-\eps \leq f(x)\leq 1\,, && 1-\eps \leq h(x,\nu)\leq 1 && \forall\ |x|\geq R-2,\ \forall\,\nu\in\S^{n-1}\,.
\end{align}
As an obvious consequence of Definition~\ref{meandensity}, if a set $F$ has weighted volume equal to $\omega_n$, then its $(f,h)$-mean density is smaller than $1$ if and only if its weighted perimeter is less than $n\omega_n$. Since by definition
\begin{align*}
|B|_f = \omega_n - |B|_{1-f}\leq \omega_n\,, && P_h(B) = n\omega_n - P_{1-h}(B)\leq n\omega_n\,,
\end{align*}
the proof is trivially concluded with $F=B$ if $|B|_f=\omega_n$. As a consequence, we can assume without loss of generality that $|B|_f<\omega_n$, or equivalently that $|B|_{1-f}>0$. Our strategy will be to obtain $F$ slightly enlarging $B$, in order to get a set whose volume is exactly $\omega_n$, but in such a way that its perimeter remains smaller than $n\omega_n$.

\step{1}{The case of a radial density $h$.}
Let us start assuming that the density $h$ is radial, that is, for every $x\in\R^n,\,\nu\in\S^{n-1}$ and for every rotation $\rho:\R^n\to\R^n$ centred at the origin, one has $h(x,\nu)=h(\rho(x),\rho(\nu))$. Let us select a hyperplane $H$ passing through the origin and the centre $R\theta$ of the ball $B$, let us call $H^\pm$ the two half-spaces in which $\R^n$ is divided by $H$, and let $B_\pm$ and $\partial_\pm B$ be the two halves of $B$ and of $\partial B$ according to~(\ref{eq:udlrset}). Let us now consider the circle $\SS_1\subseteq \S^{n-1}$ which contains the direction $\theta$ and the direction orthogonal to $H$. For every small angle $\sigma\in \SS_1$, let us call $\rho_\sigma:\R^n\to\R^n$ the rotation of angle $\sigma$ centred at the origin, and for every small $\delta>0$, let us call $F_\delta$ the set given by
\[
F_\delta := \bigcup_{0<\sigma<\delta} \rho_\sigma(B)\,.
\]
Thanks to~(\ref{farenough}), a simple integration ensures that
\[
|F_\delta|_f - |B|_f \geq \omega_{n-1}(R-1)(1-\eps)\delta\,.
\]
Since again~(\ref{farenough}) implies that $|B|_{1-f}\leq \eps\omega_n$, and since $R\gg 1$, by continuity there is some $\bar\delta\ll \eps$ such that $|F_{\bar\delta}|_f=\omega_n$. In particular, the above estimate ensures that
\begin{equation}\label{estihatdelta}
\bar\delta \leq \frac{(1-\eps)^{-1}|B|_{1-f}}{\omega_{n-1}(R-1)}\,.
\end{equation}
We set then $F:=F_{\bar\delta}$, and to conclude the proof of this step we only have to check that $P_h(F)\leq n\omega_n$. Let us define the $(n-1)$-dimensional set $\Gamma$ so that
\[
\partial F = \partial_- B \cup \rho_{\bar\delta}(\partial_+ B) \cup \Gamma\,.
\]
For every $x\in\partial B$ and every $y\in\partial F$ let us denote by $\nu(x)$ and $\nu(y)$ the normal vectors at $x$ and $y$, with respect to $\partial B$ and to $\partial F$ respectively. The radial assumption on $h$ implies that
\[
\int_{\partial_+ B} h(x,\nu(x))\, d\H^{n-1}(x) = \int_{\rho_{\bar\delta}(\partial_+ B)} h(y,\nu(y))\, d\H^{n-1}(y)\,,
\]
so that
\begin{equation}\label{uli}
P_h(F) =P_h(B)+P_h(\Gamma)\,.
\end{equation}
However, since $h\leq 1$, a trivial geometric argument, together with~(\ref{estihatdelta}) and~(\ref{eq:pallalontanobelow}), ensures that
\begin{equation}\label{exasbe}\begin{split}
\H^{n-1}_h(\Gamma) &\leq \H^{n-1}(\Gamma)
\leq (R+1) (n-1)\omega_{n-1} \bar\delta
\leq \frac{R+1}{R-1}\, (n-1)(1-\eps)^{-1}|B|_{1-f}\\
&\leq \frac{R+1}{R-1}\, \frac{(n-1)(1-\eps)^{-1}}{n-1+2\eps n}\, P_{1-h}(B)
\leq P_{1-h}(B)\,,
\end{split}\end{equation}
where the last inequality holds true because $\eps\ll 1$ and $R\gg n/\eps$. Together with~(\ref{uli}), this implies $P_h(F)\leq P_h(B)+P_{1-h}(B)=P_1(B)=n\omega_n$ and, as noticed before, this gives the thesis.

\step{2}{The general case.}
We now start showing the thesis removing the assumption that $h$ is radial. We introduce the radial averages $f_r$ and $h_r$ of $f$ and $h$ as follows,
\begin{align*}
f_r(x) = \intmed_{\alpha\in\S^{n-1}} f(\rho_\alpha(x))\,d\H^{n-1}\,, &&
h_r(x,\nu) = \intmed_{\alpha\in\S^{n-1}} h(\rho_\alpha(x),\rho_\alpha(\nu))\,d\H^{n-1}\,,
\end{align*}
where for every $\alpha\in\S^{n-1}$ we denote again by $\rho_\alpha$ the rotation of angle $\alpha$ centred at the origin. Notice that the densities $f_r$ and $h_r$ are radial by construction, and they also converge to $1$ from below and satisfy~(\ref{eq:HPbelow}), since so do $f$ and $h$. As a consequence, we can apply Lemma~\ref{prop:fhbelowstime} to get a ball $B$, arbitrarily far from the origin, such that
\begin{equation}\label{eq:plrr}
P_{1-h_r}(B) \geq \big(n-1+2\eps n \big) |B|_{1-f_r}\,.
\end{equation}
Let us call $R\gg 1$ the distance of the centre of $B$ from the origin; notice that, since $h_r$ is radial, for every $\theta\in\S^{n-1}$ the ball $B^{R\theta}_1$ satisfies $P_{1-h_r}(B^{R\theta}_1)=P_{1-h_r}(B)$. We claim that, for every $1\leq k\leq n-1$, there exists a $k$-dimensional sphere $\SS_k\approx \S^k$ contained in $\S^{n-1}$ such that
\begin{equation}\label{assstep2}
\intmed_{\SS_k} P_{1-h}(B_1^{R\theta}) - (n-1+2\eps n) |B_1^{R\theta}|_{1-f}\, d\H^k(\theta) \geq 0\,.
\end{equation}
We can easily prove this claim by induction. In fact, for $k=n-1$ it is an obvious consequence of~(\ref{eq:plrr}) with $\SS_k=\S^{n-1}$, since
\begin{align*}
\intmed_{\S^{n-1}} P_{1-h}(B_1^{R\theta})\, d\theta = P_{1-h_r}(B)\,, &&
\intmed_{\S^{n-1}} |B_1^{R\theta}|_{1-f}\, d\theta= |B|_{1-f_r}\,.
\end{align*}
Moreover, assume the property to be true for some $k\geq 2$ and some $k$-dimensional sphere $\SS_k\subseteq\S^{n-1}$. Then, for every $\theta\in \SS_k$ let us call $\SS(\theta)\approx \S^{k-1}$ the $(k-1)$-dimensional sphere contained in $\SS_k$ orthogonal to $\theta$. By homogeneity, we have
\begin{align*}
\intmed_{\SS_k} \bigg( \intmed_{\SS(\theta)} P_{1-h}(B_1^{R\alpha})\, d\H^{k-1}(\alpha)\bigg)\, d\H^k(\theta)&= \intmed_{\SS_k} P_{1-h}(B_1^{R\theta})\,d\H^k(\theta)\,, \\
\intmed_{\SS_k} \bigg( \intmed_{\SS(\theta)}  |B_1^{R\alpha}|_{1-f}\, d\H^{k-1}(\alpha)\bigg)\, d\H^k(\theta)&= \intmed_{\SS_k} |B_1^{R\theta}|_{1-f}\,d\H^k(\theta)\,,
\end{align*}
so the validity of~(\ref{assstep2}) for dimension $k$ and the sphere $\SS_k$ immediately implies the existence of some $\theta\in \SS_k$ such that~(\ref{assstep2}) holds also for dimension $k-1$ and the sphere $\SS_{k-1}=\SS(\theta)$. The claim is then established, in particular we will use a circle $\SS_1\approx \S^1$ for which~(\ref{assstep2}) holds true with $k=1$. As in Step~1, we fix $\theta\in \SS_1$ and, for every small $\delta>0$, we define
\begin{equation}\label{defetheta}
F^\theta_\delta := \bigcup_{0<\sigma<\delta} \rho_\sigma(B^{R\theta}_1)\,.
\end{equation}
Since we can assume that $1-\eps\leq f,\,h\leq 1$ around $B$ up to have taken $R$ large enough, we get the existence of a unique $\bar\delta=\bar\delta(\theta)$ such that $|F(\theta)|_f=\omega_n$, where for simplicity of notation we write $F(\theta):=F^\theta_{\bar\delta(\theta)}$; in addition, the analogous of estimate~(\ref{estihatdelta}) holds, that is,
\begin{equation}\label{newestihatdelta}
\bar\delta(\theta) \leq \frac{(1-\eps)^{-1}|B^{R\theta}_1|_{1-f}}{\omega_{n-1}(R-1)}\,.
\end{equation}
Let us now define the map $\tau:\SS_1\to \SS_1$ by $\tau(\theta)=\theta+\bar\delta(\theta)$, which is by construction a strictly increasing bijection of $\SS_1$ onto itself (keep in mind that $f\geq 1-\eps$, so in particular $f$ is strictly positive!).\par

Observe that, to conclude the thesis, it is enough to find some $\theta\in \SS_1$ such that the set $F(\theta)$ has $h$-perimeter at most $n\omega_n$, since each $F(\theta)$ has $f$-volume equal to $\omega_n$. In particular, if there exists $\theta\in \SS_1$ such that $\bar\delta(\theta)=0$, then the set $F(\theta)$ is exactly a ball of radius $1$, so its $h$-perimeter is clearly at most $n\omega_n$, being $h\leq 1$, and the proof is already concluded. We can then assume without loss of generality that $\bar\delta(\theta)>0$ for every $\theta$.\par

With a small abuse of notation, for every $\alpha\in \SS_1$ let us call $\partial_- B^\alpha$ the ``lower'' half-sphere of $\partial B^{R\alpha}_1$, that is, the set of the points in $\partial B^{R\alpha}_1$ whose direction, once projected onto the $2$-dimensional plane containing $\SS_1$, is smaller than $\alpha$ (this makes sense since $R\gg 1$, so the directions of all the projections of the points in $\partial B^{R\alpha}_1$ are almost exactly $\alpha$). Similarly, we call $\partial_+ B^\alpha$ the ``upper'' half-sphere of $\partial B^{R\alpha}_1$. Let us now fix some $\theta\in \SS_1$, and let $\zeta\ll \bar\delta(\theta)$. We can define the sets
\begin{align*}
A:=F(\theta)\setminus F(\theta+\zeta) =\bigcup\limits_{0<\sigma<\zeta} \partial_- B^{\theta+\sigma}\,, && 
C:= F(\theta+\zeta)\setminus F(\theta)=\bigcup\limits_{0<\sigma<\tau(\theta+\zeta)-\tau(\theta)} \partial_+ B^{\theta+\bar\delta(\theta)+\sigma}\,.
\end{align*}
Since $|F(\theta)|_f=|F(\theta+\zeta)|_f$, we deduce $|A|_f=|C|_f$. On the other hand, by immediate geometrical arguments, for the Euclidean volume we have
\[
\frac{|C|_1}{|A|_1} = \frac{\tau(\theta+\zeta) -\tau(\theta)}{\zeta}\,,
\]
and since $1-\eps\leq f\leq 1$ we deduce
\[
1-\eps \leq \frac{\tau(\theta+\zeta)-\tau(\theta)}\zeta\leq\frac 1{1-\eps} \,.
\]
This implies that the function $\tau: \SS_1\to \SS_1$ is bi-Lipschitz, with $1-\eps\leq\tau'\leq (1-\eps)^{-1}$. A simple change of variables yields
\[
\intmed_{\theta\in \SS_1} \H^{n-1}_{1-h}(\partial_+ B^\theta)\, d\theta 
=\intmed_{\alpha\in \SS_1} \H^{n-1}_{1-h}(\partial_+ B^{\tau(\alpha)}) \tau'(\alpha)\, d\alpha
\leq (1-\eps)^{-1}\intmed_{\alpha\in \SS_1} \H^{n-1}_{1-h}(\partial_+ B^{\tau(\alpha)}) d\alpha\,.
\]
Keeping in mind~\eqref{assstep2} we get then
\[\begin{split}
0 &\leq \intmed_{\SS_1} P_{1-h}(B_1^{R\theta}) - (n-1+2\eps n) |B_1^{R\theta}|_{1-f}\, d\theta \\
&=\intmed_{\SS_1}\H^{n-1}_{1-h}(\partial_+B^\theta)+\H^{n-1}_{1-h}(\partial_-B^\theta)-\big(n-1+2\eps n\big)|B^{R\theta}_1|_{1-f} \,d\theta\\
&\leq \intmed_{\SS_1} (1-\eps)^{-1}\H^{n-1}_{1-h}(\partial_+B^{\tau(\theta)})+\H^{n-1}_{1-h}(\partial_-B^\theta)-\big(n-1+2\eps n\big)|B^{R\theta}_1|_{1-f}\,d\theta\,,
\end{split}\]
hence we find some $\bar\theta\in\S^1$ such that
\begin{equation}\label{fff}
\H^{n-1}_{1-h}(\partial_+B^{\tau(\bar\theta)}) + \H^{n-1}_{1-h}(\partial_-B^{\bar\theta}) \geq (1-\eps)\big(n-1+2\eps n\big)|B^{R\bar\theta}_1|_{1-f}\,.
\end{equation}
Let us then call $F=F(\bar\theta)$, and let us write $\partial F = \partial_- B^{\bar\theta}\cup \partial^+ B^{\tau(\bar\theta)}\cup\Gamma$. Keeping in mind~(\ref{newestihatdelta}), arguing exactly as in~(\ref{exasbe}) we get
\[
\H^{n-1}_h(\Gamma)\leq \frac{R+1}{R-1}\, (n-1)(1-\eps)^{-1}|B^{R\bar\theta}_1|_{1-f}\,,
\]
which by~(\ref{fff}) gives
\[\begin{split}
P_h(F)&= \H^{n-1}_h(\partial_- B^{\bar\theta}) + \H^{n-1}_h (\partial_+ B^{\tau(\bar\theta)}) + \H^{n-1}_h(\Gamma)\\
&= n\omega_n -\H^{n-1}_{1-h}(\partial_- B^{\bar\theta}) - \H^{n-1}_{1-h} (\partial_+ B^{\tau(\bar\theta)}) + \H^{n-1}_h(\Gamma)\\
&\leq n\omega_n -(1-\eps)\big(n-1+2\eps n\big)|B^{R\bar\theta}_1|_{1-f}+\frac{R+1}{R-1}\, (n-1)(1-\eps)^{-1}|B^{R\bar\theta}_1|_{1-f}
\leq n\omega_n\,,
\end{split}\]
where the last inequality holds true up to have taken $\eps\ll 1/n$ and $R\gg n/\eps$. As noticed above, since $|F|_f=\omega_n$, this concludes the thesis.
\end{proof}

\subsection{The case of densities converging from above\label{secabo}}

Let us start this section by showing the existence of a ball, far away from the origin, for which the perimeter can be controlled by the volume, where both perimeter and volume are intended with respect to the density $\tilde h$ defined in~(\ref{deft}).

\begin{prop}\label{prop:disuguaglianzaabove}
Let $h$ be a function converging to $1$ such that $h^+$ converges from above to $1$ and $\int_R^{\infty }\tilde h_r(t)\,dt = +\infty$ for every $R\gg 1$. Then, for every $\eps>0$ there exists a ball $B$ of radius $1$ arbitrarily far from the origin such that
\begin{equation}\label{eq:plo}
P_{\tilde h}(B) \leq (n+\eps)|B|_{\tilde h}\,.
\end{equation}
Moreover the same holds for the radial average $\tilde h_r$ of $\tilde h$.
\end{prop}
\begin{proof}
\step{1}{Reduction to the radial case.}
Assume that the claim holds for radial densities. Recall the radial average $\tilde h_r$ of $\tilde h$ is
\[
\tilde h_r(t) = \intmed_{\partial B_{|t|}} \tilde h(x)\, d \H^{n-1}(x)\,.
\]
By assumption, for any positive $\eps$ we can find a ball $B$ arbitrarily far from the origin such that $P_{\tilde h_r}(B) \le (n+\eps)|B|_{\tilde h_r}$. Let us call for a moment $B^\theta$ the ball obtained from $B$ after a rotation of an angle $\theta\in\S^{n-1}$ around the origin. All these balls share the same perimeter and volume with respect to the radial density $\tilde h_r$, but not with respect to $\tilde h$. Since by definition
\begin{align*}
P_{\tilde h_r} (B) = \intmed_{\S^{n-1}} P_{\tilde h}(B^\theta) \,d \H^{n-1}(\theta)\,, && |B|_{\tilde h_r} = \intmed_{\S^{n-1}} |B^\theta|_{\tilde h} \,d \H^{n-1}(\theta)\,,
\end{align*}
there obviously exists some angle $\bar \theta$ for which $P_{\tilde h}(B^{\bar\theta}) \leq (n+\eps)|B^{\bar\theta}|_{\tilde h}$.
\step{2}{Proof of the radial case.}
Thanks to Step~1, we can assume without loss of generality that $\tilde h$ is radial. As a consequence, we can write the perimeter and the volume of a ball $B_R$ centred at distance $R$ from the origin as
\begin{align}\label{eq:radialPV}
P_{\tilde h}(B_R) = \int_{-1}^1 \alpha_R(t) \tilde h(t+R)\, dt\,, && |B_R|_{\tilde h}=\int_{-1}^1 \beta_R(t) \tilde h(t+R)\, dt\,,
\end{align}
where the exact value of $\alpha_R$ and of $\beta_R$ can be computed. However, for our purposes it is sufficient to observe that $\alpha_R$ and $\beta_R$ uniformly converge, as $R$ goes to infinity, to the functions $\alpha$ and $\beta$ corresponding to the flat layers, that is,
\begin{align*}
\alpha(t) = (n-1) \omega_{n-1} (1-t^2)^{\frac{n-3}2}\,, && \beta(t) = \omega_{n-1} (1-t^2)^{\frac{n-1}2}\,.
\end{align*}
Notice that $\int_{-1}^1 \alpha(t) -n\beta(t)\, dt =0$.
Let us argue now by contradiction and suppose that for every $R$ greater or equal than some $\overline R$ inequality~\eqref{eq:plo} is false for the set $B_R$. Integrating the opposite inequality over $(R_1,R_2)$ with $\overline R\ll R_1 \ll R_2$ we get
\[
\int_{R_1}^{R_2} P_{\tilde h}(B_R)\, dR > (n+\eps) \int_{R_1}^{R_2} |B_R|_{\tilde h}\,,
\]
which by~\eqref{eq:radialPV} can be written as
\[
\int_{R_1}^{R_2} \int_{-1}^1 \alpha_R(t) \tilde h(t+R)\, dt\, dR > (n+\eps) \int_{R_1}^{R_2} \int_{-1}^1 \beta_R(t) \tilde h(t+R)\, dt\, dR\,.
\]
Since $\alpha_R /\alpha \to 1$ and $\beta/\beta_R \to 1$ uniformly as $R\to\infty$, up to have taken $\overline R$ big enough we deduce
\begin{equation}\label{tph}
\int_{R_1}^{R_2} \int_{-1}^1 \alpha(t) \tilde h(t+R)\, dt\, dR > \left(n+\frac{\eps}{2}\right) \int_{R_1}^{R_2} \int_{-1}^1 \beta(t) \tilde h(t+R)\, dt\, dR\,.
\end{equation}
Notice that the left-hand side equals
\[
\int_{R_1-1}^{R_1+1} \tilde h(R) \int_{-1}^{R-R_1} \alpha(t)\, dt\, dR +\int_{R_1+1}^{R_2 -1}\tilde h(R)\int_{-1}^1 \alpha(t)\, dt\, dR +\int_{R_2-1}^{R_2+1} \tilde h(R) \int_{R-R_2}^1 \alpha(t)\, dt\, dR\,,
\]
which is clearly smaller than
\[
K + \int_{R_1+1}^{R_2 -1}\tilde h(R)\int_{-1}^1 \alpha(t)\, dt\, dR\,,
\]
where $K$ is a constant not depending on $R_1$ and $R_2$. Similarly, the right-hand side of~(\ref{tph}) equals
\[
\int_{R_1-1}^{R_1+1} \tilde h(R) \int_{-1}^{R-R_1} \beta(t)\, dt\, dR +\int_{R_1+1}^{R_2 -1}\tilde h(R)\int_{-1}^1 \beta(t)\, dt\, dR +\int_{R_2-1}^{R_2+1} \tilde h(R) \int_{R-R_2}^1 \beta(t)\, dt\, dR\,,
\]
hence it is bigger than
\[
-K + \int_{R_1+1}^{R_2 -1}\tilde h(R)\int_{-1}^1 \beta(t)\, dt\, dR\,.
\]
As a consequence, from~(\ref{tph}) we derive that, for every $\overline R\ll R_1\ll R_2$ one has
\[
K + \int_{R_1+1}^{R_2 -1}\tilde h(R)\int_{-1}^1 \alpha(t)\, dt\, dR > \bigg(n+\frac \eps 2\bigg) \bigg(-K + \int_{R_1+1}^{R_2 -1}\tilde h(R)\int_{-1}^1 \beta(t)\, dt\, dR\bigg)\,,
\]
hence
\[
\int_{R_1+1}^{R_2 -1}\tilde h(R)\int_{-1}^1 \alpha(t) - \bigg(n+\frac \eps 2\bigg) \beta (t)\, dt\, dR > - \widetilde K
\]
for some constant $\widetilde K$, again not depending on $R_1$ or $R_2$. Keeping in mind that, as noticed above, $\int_{-1}^1 \alpha - n\beta =0$, the last inequality reduces itself to
\[
- \frac \eps 2\, \int_{R_1+1}^{R_2 -1}\tilde h(R)\int_{-1}^1  \beta (t)\, dt\, dR > - \widetilde K\,.
\]
And finally, since $\int_{-1}^1 \beta>0$, up to fixing a large $R_1$ and then a much larger $R_2$, we get a contradiction with the assumption that $\int_R^{+\infty} \tilde h=+\infty$. The thesis is then proved.
\end{proof}

\begin{lem}\label{lem:fhabovestime}
Let $h$ converge to $1$, and let $f$ and $h$ converge both from above to $1$ and satisfy~\eqref{eq:HPabove}, and $\int_R^{\infty }\tilde h_r(t)\,dt = +\infty$ for every $R\gg 1$. Then, there exists a small constant $\eps>0$ such that there is a ball $B$ of radius $1$ arbitrarily far from the origin with
\begin{equation}\label{eq:pallalontanoabove}
P_{\tilde h}(B) \leq \big(n-1-\eps\big) |B|_{\tilde f}\,.
\end{equation} 
Moreover, the same holds for the radial averages $\tilde h_r$ and $\tilde f_r$ in place of $\tilde h$ and $\tilde f$.
\end{lem}
\begin{proof}
By~\eqref{eq:HPabove}, there exist a small positive constant $\delta$ and a large constant $\overline R$ such that, for every $x\in\R^n$ with $|x|>\overline R$, one has
\begin{equation}\label{average}
\tilde f(x) \geq \frac{n}{n-1} \, (1+\delta)^2\tilde h\,.
\end{equation}
Applying Proposition~\ref{prop:disuguaglianzaabove}, we get a ball $B$ with radius $1$ and distance from the origin larger than $\overline R$ such that $P_{\tilde h}(B)\leq n |B|_{\tilde h}(1+\delta)$. As a consequence, for the same ball $B$ we obtain
\[
P_{\tilde h} (B) \leq (n-1) (1+\delta)^{-1} |B|_{\tilde f} \leq \big(n-1-\eps\big) |B|_{\tilde f} \,,
\]
where the last inequality is true for a suitably small $\eps>0$, depending on $\delta$, hence ultimately on the inequality~(\ref{eq:HPabove}). Notice that the last inequality coincides with~(\ref{eq:pallalontanoabove}), which is then obtained.\par
The very same proof works with the radial averages $\tilde h_r$ and $\tilde f_r$ in place of $\tilde h$ and $\tilde f$, since Proposition~\ref{prop:disuguaglianzaabove} is stated also with $\tilde h_r$ in place of $\tilde h$, and the inequality~(\ref{average}) clearly holds true also with $\tilde f_r$ and $\tilde h_r$ in place of $\tilde f$ and $\tilde h$, as one obtains simply by an integration of~(\ref{average}) over $\S^{n-1}$.
\end{proof}

We can now conclude this section by giving the proof of Proposition~\ref{prop:SetMD<1above}.

\begin{proof}[Proof of Proposition~\ref{prop:SetMD<1above}]
First of all, as in the proof of Proposition~\ref{prop:SetMD<1below}, up to a homothety we can assume that $m=\omega_n$. In other words, we are looking for a set $F$ arbitrarily far from the origin with weighted volume $\omega_n$, and weighted perimeter smaller or equal than $n\omega_n$. Let us apply Lemma~\ref{lem:fhabovestime}, so to get a constant $\eps>0$. Let now $\eta$ be a small positive constant, depending on $\eps$ and to be specified later. Since $f$ and $h$ are converging to $1$, there is $\overline R\gg 1$ such that for every $x\in\R^n$ with $|x|>\overline R$ and every $\nu\in\S^{n-1}$ one has $1-\eta\leq h(x,\nu)\leq 1+\eta$ and $1\leq f(x)\leq 1+\eta$; in fact, keep in mind that $f\geq 1$ because $f$ is converging to $1$ from above, while it is not necessarily $h\geq 1$ since $h^+$, and not $h$, is converging to $1$ from above. Lemma~\ref{lem:fhabovestime} applied to $\tilde h_r$ and $\tilde f_r$ provides us with a ball $B=B^{R\theta}$, with radius $1$ and centred at $R\theta$, such that $R>\overline R+1$ and
\begin{equation}\label{dontfor}
P_{\tilde h_r}(B) \leq (n-1-\eps) |B|_{\tilde f_r}\,.
\end{equation}
Our aim is to use the ball $B$ to find the searched set $F$. The proof is split in three steps.
\begin{figure}[t]
\centering
\begin{tikzpicture}[line cap=round,line join=round,>=triangle 45,x=1.5cm,y=1.5cm]
\clip(-1.,-0.5) rectangle (4.7,2.5);
\draw [shift={(3.717149082283583,1.477431115171978)},line width=0.5pt,color=cqcqcq,fill=cqcqcq,fill opacity=1.0]  (0,0) --  plot[domain=-2.4850333052337543:0.02206894073719071,variable=\t]({1.*0.5*cos(\t r)+0.*0.5*sin(\t r)},{0.*0.5*cos(\t r)+1.*0.5*sin(\t r)}) -- cycle ;
\draw [shift={(3.8209787912652553,1.183267119758303)},line width=0.5pt,color=cqcqcq,fill=cqcqcq,fill opacity=1.0]  (0,0) --  plot[domain=0.656559348356036:3.1636615943269852,variable=\t]({1.*0.5*cos(\t r)+0.*0.5*sin(\t r)},{0.*0.5*cos(\t r)+1.*0.5*sin(\t r)}) -- cycle ;
\draw [line width=.5pt] (3.8209787912652553,1.183267119758303) circle (0.75cm);
\draw [line width=.5pt] (3.717149082283583,1.477431115171978) circle (0.75cm);
\draw [shift={(0.,0.)},line width=0.15pt]  (0,0) --  plot[domain=0.3003104534660773:0.3783178356271492,variable=\t]({1.*4.*cos(\t r)+0.*4.*sin(\t r)},{0.*4.*cos(\t r)+1.*4.*sin(\t r)}) -- cycle ;
\draw (0.8,0.7) node[anchor=north west] {$\bar\delta$};
\draw (3.7,1.6) node[anchor=north west] {$F$};
\draw[line width=0.1pt,->](0.6,0)--(1.4,0.25);
\draw(1,0) node[anchor=west] {$\theta$};
\draw[line width=0.1pt,->](-.3,0)--(-0.55,0.8);
\draw(-.45,.45) node[anchor=west] {$\nu$};
\draw (0.9,0.28) arc (17:21:0.94);
\end{tikzpicture}
\caption{The set $F$ of Step~1 in Proposition~\ref{prop:SetMD<1above}.}
\label{Shrink}
\end{figure}
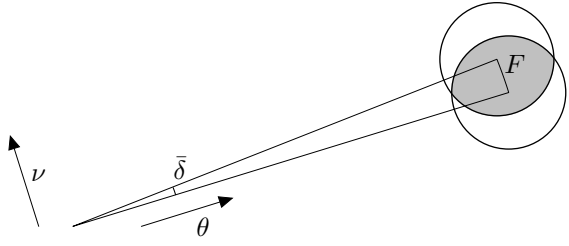
\step{1}{The case of radial densities $h$ and $f$.}
We start with the assumption that the densities $h$ and $f$ are radial, hence $\tilde h=\tilde h_r$ and $\tilde f=\tilde f_r$. Let us select a hyperplane $H$ passing through the origin and containing the centre $R\theta$ of the ball $B$, and let $\nu\in\S^{n-1}$ be the direction orthogonal to $H$, so in particular $\nu$ is orthogonal to $\theta$. With a small abuse of notation, for every $\delta\in\R$ let us now denote by ``$\theta+\delta$'' the angle obtained by $\theta$ after a rotation of angle $\delta$ in the circle contained in $\S^{n-1}$ and containing $\theta$ and $\nu$: formally speaking, we denote by $\theta+\delta$ the direction $\theta \cos \delta +\nu \sin\delta\in\S^{n-1}$. Writing for simplicity $B$ in place of $B^{R\theta}$ and $B_\delta$ in place of $B^{R(\theta+\delta)}$, we observe that by continuity there exists a small positive number $\bar\delta$ such that the set $F= B\cap B_{\bar\delta}$ satisfies $|F|_f=\omega_n$ (see Fig.~\ref{Shrink}). Keeping in mind that $f\leq 1+\eta$ on $B$, by an immediate geometric argument we have
\[
|B|_{\tilde f} = |B|_f - \omega_n = |B|_f - |F|_f \leq (1+\eta) \omega_{n-1} (R+1) \bar{\delta}\,,
\]
which implies the following lower bound for $\bar \delta$,
\begin{equation}\label{lbbd}
\bar \delta \geq \frac{1- 2\eta}{\omega_{n-1}(R+1)} \,|B|_{\tilde f}\,.
\end{equation}
Let us now call $H^\pm$ the two half-spaces in which $\R^n$ is divided by $H$, the half-space $H^+$ being the one which contains $R(\theta+\delta)$ for small positive $\delta$. We can write $\partial B=\partial^+ B \cup \partial^- B$, where $\partial^\pm B=\partial B\cap H^\pm$. Similarly, let us call $H_{\bar\delta}$ the hyperplane passing through the origin and orthogonal to the direction $\theta+(\bar\delta+\pi/2)$, and let $H_{\bar\delta}^\pm$ be the two half-spaces in which $H_{\bar\delta}$ divides $\R^n$, being $H_{\bar\delta}^+$ the one containing $R(\theta+\delta)$ for $\delta$ slightly bigger than $\bar\delta$. Notice that $H_{\bar\delta}$ passes through the centre of the ball $B_{\bar\delta}$, hence we can also split $\partial B_{\bar\delta}=\partial^+ B_{\bar\delta} \cup \partial^- B_{\bar\delta}$, where $\partial^\pm B_{\bar\delta} = \partial B_{\bar\delta} \cap H_{\bar\delta}^\pm$.\par

Let us now notice that $\partial F=\partial (B\cap B_{\bar\delta})=(\partial B\cap B_{\bar\delta}) \cup (B\cap \partial B_{\bar\delta})\subseteq \partial^+ B \cup \partial^- B_{\bar\delta}$. More precisely, we can write $\partial F =\partial^+ B \cup \partial^- B_{\bar\delta} \setminus \Gamma$, and a simple geometric consideration, also by~(\ref{lbbd}), ensures that
\[
\H^{n-1}(\Gamma) \geq (n-1)\omega_{n-1} \bar \delta (R-1)
\geq (n-1)(R-1) \,\frac{1- 2\eta}{(R+1)} \,|B|_{\tilde f}\,.
\]
Since $h$ is radial, we have $\H^{n-1}_h(\partial^- B)=\H^{n-1}_h(\partial^- B_{\bar\delta})$. Therefore, being $h\geq 1-\eta$ on $\partial F$ and by~(\ref{dontfor}) we can estimate
\begin{align*}
P_h(F) &= \H^{n-1}_h(\partial^+ B) + \H^{n-1}_h(\partial^- B_{\bar\delta}) - \H^{n-1}_h(\Gamma) = P_h(B) - \H^{n-1}_h(\Gamma)
\leq P_{1+\tilde h}(B) - \H^{n-1}_h(\Gamma)\\
&\leq n\omega_n + P_{\tilde h}(B) -(1-\eta)\,\frac{R-1}{R+1}(n-1)(1-2\eta)|B|_{\tilde f} \\
&\leq n\omega_n + \left (n-1-\eps -\frac{R-1}{R+1}(n-1)(1-3\eta)\right)|B|_{\tilde f} \leq n\omega_n\,,
\end{align*}
where the last inequality is true up to have chosen first $\eta$ small enough with respect to $\eps$, and then $R$ big enough, also recalling that $|B|_{\tilde f}\geq 0$. The set $F$ is then the desired set, and the proof in this case is concluded.

\step{2}{The general case in dimension $2$.} 
Let us now assume that $n=2$. Then, as in Step~1, we notice that for every $\theta\in\S^1$ the function $\delta\mapsto |B^{R\theta}\cap B^{R(\theta+\delta)}|_f$ is strictly decreasing for small positive $\delta$, and in particular there is a unique $\bar\delta=\bar\delta(\theta)$ such that the set $F^\theta=B^{R\theta}\cap B^{R(\theta+\bar\delta(\theta))}$ satisfies $|F^\theta|_f=\omega_n$. Moreover, exactly as in Step~1 we have the analogous of the estimate~(\ref{lbbd}) for any $\bar\delta(\theta)$, with $\theta\in\S^1$, which now reads as\par
\begin{equation}\label{lbbd2}
\bar \delta(\theta) \geq \frac{1- 2\eta}{\omega_{n-1}(R+1)} \,|B^{R\theta}|_{\tilde f}\,.
\end{equation}
We can assume without loss of generality that $\bar\delta(\theta)>0$ for every $\theta$. Indeed, if for some $\theta$ we have $\bar\delta(\theta)=0$, it means that $|B^{R\theta}|_f=\omega_n = |B^{R\theta}|_1$. Since $f\geq 1$, this means that $f\equiv 1$ on $B$, hence by lower semicontinuity and estimate~(\ref{eq:HPabove}) we get $\tilde h=0$ on $\partial B^{R\theta}$, thus $h\leq h^+=1$ on $\partial B^{R\theta}$, and finally $P_h(B^{R\theta})\leq n\omega_n$, thus the thesis is already obtained.\par

We now define the map $\tau :\S^1 \to \S^1$ as $\tau(\theta)=\theta +\bar \delta(\theta)$, and notice that by construction it is a strictly increasing bijection. Let us now fix some $\theta\in\S^1$, and any $\zeta \ll \tau(\theta)-\theta = \bar \delta(\theta)$, and call
\begin{align*}
A=  F^{\theta +\zeta} \setminus B^{R\theta}\,, && C=  F^\theta \setminus B^{R\tau(\theta+\zeta)}\,.
\end{align*}
Notice now that, by definition, $F^{\theta+\zeta}= (F^\theta\cup A) \setminus C$ and $C\subseteq F^\theta$ while $A\cap F^\theta=\emptyset$. As a consequence, since $F^\theta$ and $F^{\theta+\zeta}$ have the same weighted volume, we deduce $|A|_f=|C|_f$. On the other hand, a quick geometrical observation ensures that
\[
\frac{|C|_1}{|A|_1} = \frac{\tau(\theta+\zeta) -\tau(\theta)}{\zeta}\, (1+o(1))\,.
\]
Therefore, as $1\leq f \leq 1+\eta$ on both $A$ and $C$, we get that the map $\tau$ is bi-Lipschitz with $(1+\eta)^{-1} \leq \tau' \leq 1+\eta$. Since all the sets $F^\theta := F^\theta_{\bar \delta(\theta)}$ have weighted volume equal to $\omega_n$, to conclude we just have to find some angle $\bar \theta$ for which $P_h(F^{\bar \theta}) \leq n\omega_n$.\par

For every $\theta$, let us write $\partial B^{R\theta}=\partial^+ B_\theta\cup \partial^- B_\theta$, where a point $x\in \partial B^{R\theta}$ is said to belong to $\partial^+ B_\theta$ (resp., to $\partial^- B_\theta$) if its direction is larger (resp., smaller) than $\theta$. Notice that this makes sense since the radius of $B^{R\theta}$ equals $1$ while $R\gg 1$, hence the directions of all the points of $\partial B^{R\theta}$ are extremely close to $\theta$. A change of variables yields
\[
\intmed_{\S^1} \H^{n-1}_{\tilde h} (\partial^- B_\theta)\, d \theta = \intmed_{\S^1} \H^{n-1}_{\tilde h} (\partial^- B_{\tau(\nu)})\tau'(\nu)\, d \nu \ge (1+\eta)^{-1} \intmed_{\S^1} \H^{n-1}_{\tilde h} (\partial^- B_{\tau(\theta)})\, d \theta\,.
\]
Using~\eqref{dontfor} we find
\begin{equation}\label{staca}\begin{split}
0 &\ge P_{\tilde h_r} (B) - (n-1-\eps) |B|_{\tilde f_r} = \intmed_{\S^1} P_{\tilde h} (B^{R\theta}) - (n-1-\eps) |B^{R\theta}|_{\tilde f} \, d \theta\\
&= \intmed_{\S^1} \H^{n-1}_{\tilde h} (\partial^+ B_\theta)\, d \theta + \intmed_{\S^1} \H^{n-1}_{\tilde h} (\partial^- B_\theta)\, d \theta - (n-1-\eps) \intmed_{\S^1} |B^{R\theta}|_{\tilde f}\, d \theta \\
&\geq \intmed_{\S^1} (1+\eta)^{-1} \H_{\tilde h}^{n-1} (\partial^+ B_\theta\cup \partial^- B_{\tau(\theta)}) - (n-1-\eps) |B^{R\theta}|_{\tilde f}\, d \theta\,.
\end{split}\end{equation}
Then, there exists an angle $\bar{\theta}$ such that
\begin{equation}\label{nig1}
\H^{n-1}_{\tilde h} \left(\partial^+ B_{\bar \theta} \cup \partial^- B_{\tau(\bar \theta)}\right)\leq (1+\eta)(n-1-\eps)|B^{R\bar \theta}|_{\tilde f}\,.
\end{equation}
Let us now notice that $\partial F^{\bar\theta}=\partial^+ B_{\bar\theta} \cup \partial^- B_{\tau(\bar\theta)}\setminus \Gamma$, where $\Gamma= \partial^+ B_{\bar\theta} \cup \partial^- B_{\tau(\bar\theta)}\setminus \partial F^{\bar\theta}$. Keeping in mind that $h\geq 1-\eta$, a simple geometrical argument gives
\[
\H^{n-1}_h (\Gamma) \geq (1-\eta) \H^{n-1} (\Gamma) \geq (1-\eta) (n-1)\omega_{n-1} \bar\delta(\bar\theta) (R-1)\,.
\]
As a consequence, by~(\ref{nig1}) and~(\ref{lbbd2}) we have
\[\begin{split}
P_h(F^{\bar \theta}) &= \H^{n-1}_h \Big(\partial^+ B_{\bar \theta} \cup \partial^- B_{\tau(\bar\theta)}\Big) - \H^{n-1}_h (\Gamma)\\
&\leq n\omega_n + \H^{n-1}_{\tilde h}\Big(\partial^+ B_{\bar \theta} \cup \partial^- B_{\tau(\bar\theta)}\Big) -(1-\eta)(n-1)\omega_{n-1}\bar{\delta}(\bar \theta) (R-1) \\
&\leq n\omega_n  + (1+\eta)(n-1-\eps) |B^{R\bar \theta}|_{\tilde f} -
(1-3\eta)(n-1)\frac{R-1}{R+1} \,|B^{R\bar\theta}|_{\tilde f}\,.
\end{split}\]
And finally, we deduce again the searched inequality $P_h(F^{\bar\theta})\leq n\omega_n$ as soon as we choose first $\eta$ small enough depending on $\eps$, and then $R$ big enough. The proof is then concluded also in this case.

\step{3}{The general case.}
We conclude now the proof by considering the general case $n\geq 3$. A quick look at Step~2 and in particular at the key calculation~(\ref{staca}) ensures that, whatever $n\geq 2$ is, if there is some unit circle $\C\approx \S^1$ inside $\S^{n-1}$ such that
\begin{equation}\label{eq:targeteq1}
\intmed_{\C} P_{\tilde h} (B^{R\theta})\, d \H^1(\theta) \leq (n-1-\eps) \intmed_\C |B^{R\theta}|_{\tilde f}\, d \H^1(\theta)\,,
\end{equation}
then the very same proof as in Step~2 works (to make this check even simpler, in Step~2 we always used the generic letter $n$ even if we were assuming $n=2$). And in fact, for a generic $n$ the estimate~(\ref{dontfor}) can be written as
\begin{equation}\label{eq:targeteq2}
\intmed_{\S^{n-1}} P_{\tilde h} (B^{R\theta})\, d \H^{n-1}(\theta) \leq (n-1-\eps) \intmed_{\S^{n-1}} |B^{R\theta}|_{\tilde f}\, d \H^{n-1}(\theta)\,,
\end{equation}
which coincides with~(\ref{eq:targeteq1}) when $n=2$ with the only possible choice $\C=\S^1$. As a consequence, we can easily argue by induction. Fix $n\geq 3$ and suppose that for some $2\leq k\leq n-1$ there is a $k$-dimensional sphere $\SS_k\subseteq\S^{n-1}$ such that
\begin{equation}\label{generic}
\intmed_{\SS_k} P_{\tilde h} (B^{R\theta})\, d \H^k(\theta) \leq (n-1-\eps) \intmed_{\SS_k} |B^{R\theta}|_{\tilde f}\, d \H^k(\theta)\,.
\end{equation}
Then, for every $\theta\in\SS_k$ let us call $\SS(\theta)$ the set of those vectors in $\SS_k$ which are orthogonal to $\theta$. Notice that every $\SS(\theta)$ is a $(k-1)$-dimensional sphere contained in $\S^{n-1}$, and by the rotational invariance of the Hausdorff measure one has
\begin{align*}
\intmed_{\SS_k} P_{\tilde h} (B^{R\theta})\, d\H^k(\theta) &= \intmed_{\SS_k} \intmed_{\SS(\theta)} P_{\tilde h} (B^{R\sigma})\, d \H^{k-1}(\sigma) \,d\H^k(\theta)\,,\\
\intmed_{\SS_k} |B^{R\theta}|_{\tilde f}\, d\H^k(\theta) &= \intmed_{\SS_k} \intmed_{\SS(\theta)} |B^{R\sigma}|_{\tilde f}\, d \H^{k-1}(\sigma) \,d\H^k(\theta)\,.
\end{align*}
Thus, then there exists some $(k-1)$-dimensional sphere $\SS_{k-1}=\SS(\theta)\subseteq \SS_k$ for which
\[
\intmed_{\SS_{k-1}} P_{\tilde h} (B^{R\theta})\, d \H^{k-1}(\theta) \leq (n-1-\eps) \intmed_{\SS_{k-1}} |B^{R\theta}|_{\tilde f}\, d \H^{k-1}(\theta)
\]
holds. In other words, for every $2\leq k\leq n-1$ the existence of a $k$-dimensional sphere satisfying~(\ref{generic}) implies the existence also of a $(k-1)$-dimensional sphere satisfying~(\ref{generic}) with $k-1$ in place of $k$. Since this property is true with $k=n-1$ by~(\ref{eq:targeteq2}), by induction we obtain the same property with $k=1$ and some $1$-dimensional sphere, that is, we have found a circle $\C$ satisfying~(\ref{eq:targeteq1}). As noticed above, the existence of a circle $\C$ satisfying~(\ref{eq:targeteq1}) yields the thesis.
\end{proof}

\section{A counterexample\label{sec:counter}}

In this last section we give an example showing that, in Theorem~\mref{Main}, the assumption that $\int_R^\infty \tilde h_r=+\infty$ for the case of densities converging from above is really needed. This can be at first sight a bit surprising, also considering that the analogous property is not needed for the case of densities converging from below. Nevertheless, a careful look at the proof of the theorem allows to realise that the assumption is not just useful for the proof, but actually essential.\par

In our example we consider the case $n=2$ just for simplicity, but it appears clear that the very same construction can be done also for a generic $n\geq 3$. For every $x\in \R^2$ and $\nu\in\S^1$, we set
\begin{align*}
f(x)= 1+ 3\varphi(|x|)\,, && h(x,\nu) = 1+ \varphi(|x|)\,,
\end{align*}
where the function $\varphi:[0,+\infty)\to (0,+\infty)$ is given by
\[
\varphi(t) = M e^{-M(t-1)^+}
\]
for some $M$ large enough to be precised later. Notice that $f$ and $h^+=h$ are converging from above to $1$ and the assumption~(\ref{eq:HPabove}) is satisfied, but $\int_R^\infty \tilde h_r=\int_R^\infty \varphi<+\infty$. In particular, observe that $\varphi(t)=M$ for every $t\leq 1$, while for $t>1$ one has $\varphi(t)= M e^M e^{-Mt}$, hence $f$ and $h$ are converging to $1$ at exponential speed. We will show the following result.

\begin{lem}
There exists no isoperimetric set of volume $\pi$ for the densities $f$ and $h$ defined above.
\end{lem}
\begin{proof}
Since both $f$ and $h$ are converging to $1$, a unit ``ball at infinity'' has volume $\pi$ and perimeter $2\pi$, hence $\J(\pi)\leq 2\pi$. As a consequence, it is enough to prove that $P_h(E)>2\pi$ for every set $E\subseteq\R^2$ such that $|E|_f=\pi$. Of course, without loss of generality, we can assume that $E$ is smooth, and $\H^1(\partial E\cap \partial B)=0$, where $B$ is the unit ball centred at the origin. Let then $E$ be such a set, and let us write for brevity $|\cdot|$ and $P$, in place of $|\cdot|_1$ and $P_1$, to denote the Euclidean volume and perimeter. Since $|E|_f=|E| + 3|E|_\varphi$, one has $|E|<\pi$, hence the standard isoperimetric inequality gives
\begin{equation}\label{qs}
P(E) \geq 2 \sqrt \pi \sqrt{|E|} > 2|E|\,.
\end{equation}
We claim that
\begin{equation}\label{finfre}
P_\varphi(E) \geq 6 |E|_\varphi\,.
\end{equation}
Notice that, if this estimate holds, then by~(\ref{qs}) we get
\[
P_h(E) = P(E) + P_\varphi(E) 
> 2 |E| + 2 |E|_{3\varphi} = 2 |E|_f = 2\pi\,,
\]
hence the thesis is obtained. Therefore, to conclude the proof it is enough to establish~(\ref{finfre}). We divide the situation in two cases.
\step{I}{If $E\cap B=\emptyset$.}
Let us first assume that $E$ has empty intersection with the unit ball $B$. Hence, around $E$, the function $\varphi(t)$ coincides with $e^{-Mt}$ up to a multiplicative constant, which of course does not play any role in the proof of~(\ref{finfre}). In order to deal later with the case of non-empty intersection, we will now show an estimate stronger than~(\ref{finfre}), namely,
\begin{equation}\label{finfre2}
P_\varphi(E) \geq 12 |E|_\varphi\qquad \hbox{if } E\cap B=\emptyset\,.
\end{equation}
For every $s>0$, let us call $\tau(s)=\H^1(E \cap \partial B_s)$. We can assume without loss of generality that
\begin{equation}\label{wlog}
\tau(s) \leq 2\pi \qquad \forall \, s>0\,. 
\end{equation}
Indeed, if~(\ref{wlog}) fails, then we directly get $P_h(E)> P(E)>2\pi$, so the thesis is already obtained without any need of establishing~(\ref{finfre}) or~(\ref{finfre2}).\par

We now define $t_1 = \inf\{ t>0:\, \tau(t) > 1/12\} \in [1,+\infty]$. Observe that
\[
|E \cap B_{t_1}|_\varphi = \int_0^{t_1} \varphi(s)\tau(s)\, ds \leq \frac 1{12} \int_{\{s\leq t_1,\, \tau(s)>0\}} \varphi(s)\, ds
\leq \frac 1{24} \, \H^1_\varphi\big((\partial E) \cap B_{t_1}\big)\,.
\]
If $t_1=\infty$, this estimate reads as $P_\varphi(E) \geq 24 |E|_\varphi$, so even stronger than~(\ref{finfre2}). Hence, we can assume that $t_1<\infty$, which readily implies, since $\varphi$ is decreasing,
\begin{equation}\label{goodstart}
P_\varphi(E)-12 |E\cap B_{t_1}|_\varphi \geq \frac{\H^1_{\varphi}\big((\partial E) \cap B_{t_1}\big)}2 \geq \frac{\varphi(t_1)}{25}\,.
\end{equation}
Let us now write for brevity $t^+ = t_1 + (1200\pi)^{-1}$. Keeping in mind~(\ref{wlog}) and again the fact that $\varphi$ is decreasing, we have
\[
\big|E \cap (B_{t^+}\setminus B_{t_1})\big|_\varphi 
=\int_{t_1}^{t^+} \varphi(s) \tau(s)\, ds \leq \varphi(t_1) 2\pi (t^+-t_1) = \frac{\varphi(t_1)}{600}\,,
\]
while
\[
\big|E\setminus B_{t^+}|_\varphi \leq \varphi(t^+) |E\setminus B_{t^+}| \leq \pi \varphi(t^+) = \pi \varphi(t_1) e^{-\frac M{1200\pi}}\leq
\frac{\varphi(t_1)}{600}\,,
\]
where the last estimate clearly holds as soon as $M$ is big enough. The last two estimates, together with~(\ref{goodstart}), give
\[
P_\varphi(E)\geq 12 |E \cap B_{t_1}|_\varphi + \frac{\varphi(t_1)}{25} = 12 \bigg(|E \cap B_{t_1}|_\varphi + \frac{\varphi(t_1)}{300} \bigg) \geq 12 |E|_\varphi\,,
\]
hence~(\ref{finfre2}) is obtained and this step is concluded. Notice that in this step we did not use the assumption that $|E|_f =\pi$. Therefore, (\ref{finfre2}) holds for every set $E\subseteq\R^2\setminus B$, regardless of its volume.

\step{II}{If $E\cap B\neq \emptyset$.}
To finish the proof, let us assume now that $E$ has a non-empty intersection with the unit ball $B$, and let us call $E^+=E\setminus B$ and $E^-=E\cap B$. Since $\varphi\equiv M$ in $E^-$, we have
\[
|E^-| = \frac{|E^-|_\varphi}M \leq \frac \pi M\,,
\]
so that the standard isoperimetric inequality gives
\begin{equation}\label{g2}
P(E^-) \geq 2 \sqrt \pi \sqrt{|E^-|} \geq 2\sqrt M |E^-|\,.
\end{equation}
Let us now call $S=E \cap \partial B$. If $M$ is big enough, then $|E^-|\ll |B|$, hence we have
\[
\H^1\big((\partial E)\cap B\big) \geq \frac 34 \H^1(S)\,,
\]
which implies
\begin{equation}\label{g1}
\H^1\big((\partial E)\cap B\big) \geq \frac 17\, P(E^-)+ \frac12\,\H^1(S)\,.
\end{equation}
Keep in mind that, as noticed at the end of Step~I, the validity of~(\ref{finfre2}) is true for every set without intersection with $B$. As a consequence, by the fact that $E^+\cap B=\emptyset$, we get~(\ref{finfre2}) with $E^+$ in place of $E$. Then, recalling that $E$ is smooth and $\H^1(\partial E\cap \partial B)=0$, using the fact that $\varphi=M$ on the closure of $B$, hence also on $S$, and by~(\ref{g1}), (\ref{g2}) and~(\ref{finfre2}) for $E^+$, we get
\[\begin{split}
P_\varphi(E) &= \H^1_\varphi \big((\partial E)\cap B\big)+ \H^1_\varphi\big((\partial E)\setminus B\big)
= M \H^1\big((\partial E)\cap B\big)+ \H^1_\varphi\big((\partial E)\setminus B\big)\\
&\geq \frac 17\, M \,P(E^-)+ \frac12\,\H^1_\varphi(S) + \H^1_\varphi\big((\partial E)\setminus B\big)
\geq \frac 27 M \sqrt M |E^-| + \frac 12\, P_\varphi(E^+)\\
&\geq \frac 27 \sqrt M |E^-|_\varphi + 6 |E^+|_\varphi
\geq 6 |E|_\varphi\,,
\end{split}\]
where again the last estimate holds if $M$ is big enough. Summarising, we have established the validity of~(\ref{finfre}), so the thesis is obtained.
\end{proof}

\bibliographystyle{plain}

\bibliography{densities}

\end{document}